\magnification 1200

\def\Spec#1{\Phi _ {#1}}
\newskip\litemindent
\litemindent= 1.7 cm  
\def\Litem#1#2{\par\noindent\hangindent#1\litemindent
\hbox to #1\litemindent{\hfill\hbox to \litemindent
{#2 \hfill}}\ignorespaces}
\def\litem{\Litem1}

\def\frac #1#2{{{#1}\over{#2}}}
\def\rd{\hbox{\rm d}}
\def\C{{\bf C}}
\def\N{{\bf N}}
\def\Q{{\bf Q}}

\centerline {\bf Trivial Jensen measures without regularity}
\smallskip
\centerline{\bf J.~F.~Feinstein}
\bigskip
\noindent {\bf Abstract} In this note we construct Swiss cheeses $X$ such that $R(X)$ is
non-regular but such that $R(X)$ has no non-trivial Jensen measures. We
also construct a non-regular uniform algebra with compact,
metrizable character space such that every point of the character space is a
peak point.
\bigskip
In [Co] Cole gave a counterexample to the peak point conjecture
by constructing a non-trivial uniform algebra $A$ with compact, metrizable
character space $\Phi_A$ such that every point of $\Phi_A$ is a peak point
for $A$. This uniform algebra was obtained from an example of McKissick [M]
by a process of repeatedly adjoining square roots. Because McKissick's
algebra is regular,
so is this first example constructed by Cole (see [F2] and [Ka]). This leads to
the following question:
Let $A$ be a uniform algebra with compact, metrizable
character space $\Phi_A$ such that every point of $\Phi_A$ is a peak point
for $A$. Must $A$ be regular?

In this note we construct an example of a Swiss cheese
$X$ for which the uniform algebra $R(X)$ is non-regular, but such that $R(X)$
has no non-trivial Jensen measures. We then apply Cole's construction to this
example to produce an example of a non-regular uniform algebra with compact,
metrizable character space such that every point of the character space is a
peak point.
\medskip
We begin by recalling some standard facts about
Jensen measures.
\smallskip
\noindent
{\bf Notation} For a commutative Banach algebra $A$,
we denote by $\Spec{A}$ the character space of $A$.
Now suppose that $A$ is a uniform algebra on a compact
space $X$.
For $x\in X$, let $M_x$ and $J_x$ be the ideals of
functions in $A$ vanishing at $x$, and in a neighbourhood of
$x$, respectively.

\smallskip
\noindent
{\bf Definition 1} Let $A$ be a uniform algebra on a
compact space $X$, and let
$\phi \in \Spec{A}$. Then a {\it Jensen} measure for $\phi$ is a regular,
Borel probability measure $\mu$ on $X$ such that, for all $f \in A$,
$$\log |\phi(f)| \leq \int_X{\log |f(x)| \rd \mu(x)}$$
(where $\log(0)$ is defined to be $-\infty$).
\smallskip

Let $x\in X$. We say that $x$ is a {\sl point of continuity} for $A$
if there is no point $y$ of $X \setminus \{x\}$ satisfying $M_x\supseteq J_y$.
\smallskip
It is standard (see, for example, [G, p.33]) that every $\phi \in \Spec{A}$
has a Jensen measure supported on $X$, and that each such measure
{\it represents} $\phi$, i.e., for all $f \in A$,
$$\phi(f) = \int_X{f(x) \rd \mu(x)}.$$

It is elementary to see that if $x$ is a point of continuity for $A$ then the
only Jensen measure for the evaluation character at $x$ which is supported
on $X$ is the point mass at $x$. The converse is, in general, false as is easily
seen by considering the disc algebra: here the only representing measures for points
of the unit circle are point masses, while clearly there are no points of continuity.

In the case where $X = \Spec{A}$, the usual definition of {\it regularity} of $A$ is
that, for each closed subset $E$ of $X$ and each $x \in X\setminus E$
there is an $f\in A$ such
that $f(x) = 1$ and $f$ is identically $0$ on $E$. This is easily seen to be equivalent
to the fact that every point of $X$ is a point of continuity for $A$. The first example
of a regular uniform algebra which was non-trivial (i.e. not equal to the algebra of
all continuous functions on a compact space) was given in [M]. This example was modified by
O'Farrell in [O] to give an example of a regular uniform algebra which has a non-zero,
continuous point derivation of infinite order. When regularity fails, it usually fails
fairly drastically: see [FS] for results on non-regularity.

It is clear that if $A$ is a regular uniform algebra, then there are no non-trivial
Jensen measures supported on $X$. The examples constructed in this note show that
the converse is false: there are non-regular uniform algebras for which
there are no non-trivial Jensen measures.

There are many examples already known of Swiss
cheeses $X$ for which $R(X)$ has no non-trivial Jensen measures:
as well as
McKissick's original example of a non-trivial regular $R(X)$ ([M], but see
also [Ko]), other examples can be found in [H], [Bro, pp.193--195] and
also (in view of Theorem 3.8 of [G], as explained on page 64 of
[G]) in [W]. It is not
clear whether or not $R(X)$ is regular for these latter examples.

In this note we shall construct a compact plane set $X$ such that $R(X)$ has no
non-trivial Jensen measures, but such that $R(X)$ is definitely not regular.
\medskip
We shall make heavy use of the following version of McKissick's Lemma, given in
[Ko].
\smallskip
\noindent
{\bf Lemma 2} {\sl Let $D$ be an open disc in $\C$ and let $\epsilon > 0$.
Then there is a sequence $\Delta_k (k\in\N)$
of (pairwise disjoint) open discs with each
$\Delta_k \subseteq D$
such that the sum of the radii of the $\Delta_k$ is less than $\epsilon$
and such that, setting $U = \bigcup_{k\in\N}{\Delta_k}$, there is a sequence $f_n$
of rational functions with poles only in $U$ and such that $f_n$ converges uniformly on
$\C\setminus U$ to a function $F$ such that $F(z)=0$ for all $z \in \C \setminus D$ while
$F(z) \neq 0$ for all $z \in D\setminus U$.
}
\medskip
We also need the following elementary lemma, which follows from Cauchy's integral
formula applied to a sequence of compact sets obtained by deleting finitely many open
discs from the closed unit disc. These estimates are well-known:
for an explicit proof
of the estimate for the first derivative, see, for example, [F1].
\smallskip
\noindent
{\bf Notation} For a bounded, complex-valued function $f$
defined on a non-empty set $S$ we
shall denote by $|f|_S$ the uniform norm of $f$ on $S$, that is
$$|f|_S = \sup\{|f(x)|: x\in S\}.$$
For a compact plane set $X$, $R_0(X)$ is the set of restrictions to
$X$ of rational functions with poles off $X$. (So $R(X)$ is the uniform
closure of $R_0(X)$.)
\medskip
\noindent
{\bf Lemma 3} {\sl
Let $D_n$ be a sequence of open discs in $\C$ (not necessarily
pairwise disjoint), and set
$X = \bar{\Delta}\setminus \bigcup_{n=1}^{\infty}{D_n}$.
Suppose that $z\in X $.
Let $s_n$ denote the distance from $D_n$ to $z$ and $r_n$ the radius of $D_n$. We
also set $r_0 = 1$ and $s_0 = 1-|z|$. Suppose that $s_n > 0 $ for all $n$. Then,
for all $f \in R_0(X)$ and $k\geq 0$, we have

$$|f^{(k)}(z)| \leq k!  \sum_{j=0}^{\infty}{\frac{r_j}{{s_j}^{k+1}}} |f|_X.$$

}
\medskip
We shall also require the following famous lemma due originally to Denjoy [D],
under some additional hypotheses which were shown to be unnecessary in [C].
\smallskip
\noindent
{\bf Lemma 4} {\sl
Let $f$ be an infinitely differentiable function on an interval $I$ such that
$$\sum_{k=0}^{\infty}{\frac {1}{|f^{(k)}|_I^{1/k}}} = \infty.\eqno{(*)}$$
Suppose that there is an $x \in I$ such that $f^{(k)}(x) = 0$ for all
$k \geq 0$. Then $f$ is constantly $0$ on $I$.
}
\smallskip
Any algebra of infinitely differentiable (complex-valued) functions on
$I$ satisfying condition $(*)$ is thus a {\sl quasianalytic}
algebra, in the sense of Denjoy and Carleman.
\smallskip
We are now able to construct the desired Swiss cheese. The method for
ensuring quasianalyticity in the following theorem is based on
the method used by Brennan in [Bre]. However, we
use McKissick's Lemma and the Cauchy estimates above, and do not need to appeal
to the theory of the Bergman kernel.
\smallskip
\noindent
{\bf Theorem 5} {\sl
Set $I =[-1/2,1/2]$. There is a Swiss cheese $X$ with  $I \subseteq X$ such that
every point of $X\setminus I$ is a point of continuity for $R(X)$ but such that
every $f \in R(X)$ is infinitely differentiable on $I$ and satisfies condition
$(*)$ there.
}
\smallskip
\noindent
{\bf Proof} The Swiss cheese $X$ is constructed inductively, deleting a countable
collection of discs at each stage, and applying the above lemmas.
First we set $\delta_n = {1 \over{n+2}}$ for $n \in \N$, and we set
$$K_n=\{z \in \C: {\rm dist}(z,I) \leq \delta_k\}.$$
Let ${\cal S}_n$ be the set of all discs with rational radius which are
centred on points of
$\Q + i\Q$ and which do not meet $K_n$. Enumerate ${\cal S}_n$ as
$\{D_{n,1}, D_{n,2}, \dots\}$.
Set $X_0 = \bar\Delta$ (the closed unit disc).

Now set $\epsilon_1 = 1/4$. For each $k \in \N$, apply McKissick's lemma to
the discs $D_{1,k}$ to obtain a countable family of open discs ${\cal F}_{1,k}$
and a sequence of rational functions $f_{1,k,n}$ as in that lemma, such that
the sum of the radii of
the discs in ${\cal F}_{1,k}$ is less than ${\epsilon_1}\over{2^k}$, and, with
$U_{1,k}=\bigcup\{\Delta:\Delta \in {\cal F}_{1,k}\}$, the functions
$f_{1,k,n}$ have poles only in $U_{1,k}$ and converge to a function $F_{1,k}$
uniformly on
$\C \setminus U_{1,k}$, where $F_{1,k}$ vanishes identically on $\C \setminus D_{1,k}$
and is nowhere zero on $D_{1,k} \setminus U_{1,k}$.
Set $X_1 = X_0 \setminus \bigcup_{k\in\N}{U_{1,k}}.$ Then, for $f \in R_0(X_1)$, we
can use Lemma 3 to estimate the derivatives of $f$ on $I$: for $z \in I$ and
$n \in \N$ we have
$$|f^{(n)}(z)| \leq {{n!(1+\epsilon_1)}\over{\delta_1^{n+1}}}|f|_{X_1}.$$
Set
$$A_{1,n} = {{2n!(1+\epsilon_1)}\over{\delta_1^{n+1}}}.$$
Note that $\sum A_{1,n}^{-1/n}$
diverges. Choose $N_1$ such that $\sum_{n=1}^{N_1}{A_{1,n}^{-1/n}}\geq 1$.
We now move to the next stage of the construction.
Choose $\epsilon_2 > 0$ small enough that
$n!((1+\epsilon_1)/\delta_1^{n+1} + \epsilon_2/\delta_2^{n+1}) < A_{1,n}$ for
$1 \leq n \leq N_1$. For each $k\in \N$, apply McKissick's lemma
to the discs $D_{2,k}$
to obtain a countable family of open discs ${\cal F}_{2,k}$
and a sequence of rational functions $f_{2,k,n}$ as in that lemma, such that
the sum of the radii of
the discs in ${\cal F}_{2,k}$ is less than $\epsilon_2/2^k$, and, with
$U_{2,k}=\bigcup\{\Delta:\Delta \in {\cal F}_{2,k}\}$, the functions
$f_{2,k,n}$ have poles only in $U_{2,k}$ and converge to a function $F_{2,k}$
uniformly on
$\C \setminus U_{2,k}$, where $F_{2,k}$ vanishes identically on $\C \setminus D_{2,k}$
and is nowhere zero on $D_{2,k} \setminus U_{2,k}$.
Set $X_2=X_1\setminus \bigcup_{k\in\N}{U_{2,k}}$. Applying Lemma 3
we see that, for $f \in R_0(X_2)$, $z \in I$ and $n \in \N$,
$|f^{(n)}(z)| \leq n!((1+\epsilon_1)/\delta_1^{n+1} +
  \epsilon_2/\delta_2^{n+1}) |f|_{X_2}$.
Set $A_{2,n} = A_{1,n}$ for $1 \leq n \leq N_1$, and
$$A_{2,n} = 2n!\left({{{1+\epsilon_1}\over{\delta_1^{n+1}}} +
{{\epsilon_2}\over{\delta_2^{n+1}}}}\right)$$
for $n>N_1$.
Again we see that $\sum_{n=1}^\infty{A_{2,n}^{-1/n}} = \infty$, so choose
$N_2 > N_1$ such that $\sum_{n=N_1+1}^{N_2}{A_{2,n}^{-1/n}} \geq 1$.
Now choose $\epsilon_3 > 0 $ such that
$n!((1+\epsilon_1)/\delta_1^{n+1} + \epsilon_2/\delta_2^{n+1} +
\epsilon_3/\delta_3^{n+1})< A_{2,n}$ for
$1 \leq n \leq N_2$. We now proceed to choose families of discs and functions as before.
The inductive process is now clear, and produces a decreasing family of compact sets
$X_n$, arrays of functions $F_{n,k}$ and positive real numbers $A_{n,k}$ for
$n,k \in \N$ and a strictly increasing sequence of
positive integers $N_n$ such that $A_{n,k}=A_{j,k}$ whenever $1 \leq j \leq n$ and
$1 \leq k \leq N_n$. We may now define $A_k$ by setting $A_k = A_{n,k}$ whenever
$k \leq N_n$ (this is well defined by above).
We have $\sum_{k=1}^\infty{A_k}^{-1/k}=\infty$.
Set $X = \bigcap_{n=1}^\infty{X_n}$.
Then we claim $X$ is a Swiss cheese with the required properties. First we note that
for $f \in R_0(X)$, $z \in I$ and $k \in \N$ we have
$|f^{(k)}(z)| \leq A_k |f|_X$. It thus follows that every element of $R(X)$ is
infinitely differentiable on $I$ and satisfies condition $(*)$ there. Finally, for
any point $z \in X\setminus I$, we see that there is some $n$ with
$z \in X \setminus K_n$. Let $w \in X$ with $w \neq  z$. Then there is some disc
$D_{n,k}$ with $z \in D_{n,k}$ and $w \in X \setminus \overline{D_{n,k}}$. Then
$F_{n,k}|X$ is in $R(X)$ and is non-zero at $z$, but vanishes in a neighbourhood
of $w$. This shows that $J_w$ is not a subset of $M_z$. Thus each such point $z$ is
a point of continuity for $R(X)$, as required.
Q.E.D.
\medskip
\noindent
{\bf Corollary 6} {\sl
Let $X$ be the Swiss cheese constructed in the preceding Theorem.
Then $R(X)$ has no non-trivial
Jensen measures, but $R(X)$ is not regular.
}
\smallskip
\noindent
{\bf Proof} Certainly there are no non-trivial Jensen measures for points of
$X\setminus I$. Also, $R(X)$ is not regular, since none of the points of $I$ are
points of continuity for $R(X)$. It remains to show that there are no non-trivial
Jensen measures for points of $I$. However, from the standard theory of Jensen interior
(or fine interior), see page 319 of [GL]
or  numerous papers of Debiard and Gaveau,
it is clear that if $R(X)$
has any non-trivial Jensen measures, then
the set of points which have non-trivial Jensen measures is fairly
large:
certainly it
could not be contained in the interval $I$. The result follows. Q.E.D.
\medskip
\noindent
We now show that there are non-regular uniform algebras with
compact, metrizable character space for which
every point of the character space is a
peak point (which is equivalent to saying that
every maximal ideal has a bounded approximate identity).
We shall do this
by applying the following result to the algebra constructed
in Corollary 6. This result is a combination of results from
[Co] and [F2] (see also [Ka]), and is based on Cole's systems
of root extensions for uniform algebras.
\medskip
\noindent
\noindent {\bf Proposition 7}
{\sl Let $A_1$ be a uniform algebra on
a compact, Hausdorff space $X_1=\Phi_{A_1}$
such that the only Jensen measures for characters of $A_1$
are point masses.
Then there is a uniform algebra $A$ on a
compact, Hausdorff space
$X=\Phi_A$, a surjective continuous map $\pi$ from $X$ onto $X_1$ and a
bounded linear map $S: A \rightarrow A_1$ with the following
properties.

(a) Every maximal ideal
in $A$ has a bounded approximate identity.

(b) For every $f \in A_1$, $f \circ \pi$ is in $A$.

(c) If $x\in X_1$ and $g \in A$ with $g$ constantly equal to some complex
number $c$ on $\pi^{-1}(\{x\})$ then $(Sg)(x) = c$. In particular, for all
$f\in A_1$, $S(f \circ \pi) = f$.

(d) If $A_1$ is regular then so is $A$.

\noindent If $X_1$ is metrizable, then in addition to the above properties we
may also
insist that $X$ is metrizable.}

\smallskip
\noindent
Note that it is easy, using (c), to see that the converse to
(d) also holds: if $A$ is regular, then so is $A_1$.
\medskip
\noindent
{\bf Corollary 8}  {\sl There is a uniform algebra $A$ on a compact,
metrizable topological space
$X=\Spec{A}$
such that every point of $X$ is a peak point and such that
$A$ is not regular.}
\smallskip
\noindent
{\bf Proof} Let $X_1$ be the Swiss cheese constructed in Theorem 5, and
set $A_1=R(X_1)$ (so $A_1$ is not regular, but the only
Jensen measures for $A_1$ are point masses).
Now apply Proposition 7, and the remark following it, to
produce a non-regular, uniform algebra $A$ on a compact, metrizable
topological space $X=\Phi_A$ such that every maximal ideal of $A$ has
a bounded approximate identity. This algebra has
the required properties. Q.E.D.
\medskip
\noindent
We conclude with some open questions. An alternative way to construct
the example of Corollary 8 is to use Basener's simple construction of
non-trivial uniform algebras for which every point of the character space
is a peak point (announced in [Ba], but see pages
202-203 of [S] for the details). It is easy to see, as
with Cole's construction,
that Basener's example is non-regular whenever the original algebra $R(X)$ is
non-regular. This
leads to two related questions:
\smallskip\noindent
{\bf Question 1.} Is Basener's example ever regular?
\smallskip\noindent
{\bf Question 2.} Is Basener's example always regular when
the original algebra $R(X)$ is?
\smallskip\noindent
Our final question concerns quasianalyticity for $R(X)$. There are
several alternative definitions in print for a collection of functions
to be quasianalytic. In [Bre], Brennan gave an example of a
Swiss cheese $X$
such that, for all $p\geq2$, no non-zero element of
$R^p(X)$ vanishes almost everywhere on a subset of $X$ with positive area.
Here
$R^p(X)$ is the closure of $R(X)$
in $L^p(X)$ (area measure).

Let $X$ be a compact plane set with positive area.
We shall say that $R(X)$
is {\it quasianalytic} if the only function in $R(X)$ whose zero
set has positive area is
the constant function $0$.
\smallskip
\noindent
{\bf Question 3.} Let $X$ be a compact plane set with positive area.
Is it possible for $R(X)$ to be quasianalytic and yet to have no
non-trivial Jensen measures?
\medskip
The example constructed in Theorem 5 above shows that some elements of
quasianalyticity
may be introduced, but it is not clear whether this method can be extended to
answer Question 3. We may also ask the same question using a different form
of quasianalyticity, for example insisting that the only function in $R(X)$
which vanishes
identically on a non-empty, relatively open subset of $X$ is the zero function.
\medskip
I would like to take this opportunity to thank the Faculty members of the
Department of Mathematics of Brown University for many stimulating discussions
during my visit there in 1999 which were very beneficial for my work. I
would also like to thank the University of Nottingham for granting me sabbatical
leave in order to allow me to visit Brown University.
\bigskip
\centerline{\bf REFERENCES}
\smallskip
\litem{[Ba]}	Basener,~R.F., An example concerning peak points,
		{\it Notices Amer. Math. Soc} 18 (1971), 415-416.

\litem{[Bre]}	Brennan, James~E., Approximation in the mean and quasianalyticity,
		{\sl J. Functional Analysis}, 12 (1973) 307--320.
		
\litem{[Bro]}	Browder,~A., {\sl Introduction to Function Algebras},
		W.A. Benjamin, New York, 1969.

\litem{[Ca]}	Carleman,~T., Sur un th\'eor\`em de M. Denjoy,
		{\it C.R. Acad. Sci. Paris} 174 (1922), 121--124.

\litem{[Co]}	Cole,~B., {\sl One point parts and the peak point conjecture},
		Ph.D. dissertation, Yale Univ., 1968.

\litem{[De]}	Denjoy,~A., Sur les fonctions quasi-analytiques de variable
		r\'eelle, {\it C.R. Acad. Sci. Paris} 173 (1921), 1329--1331.

\litem{[F1]}	Point derivations and prime ideals in R(X),
		{\it Studia Mathematica}, 98 (1991) 235--246.

\litem{[F2]}	Feinstein,~J.F., A non-trivial, strongly regular uniform
		algebra, {\sl J. London Math. Soc.}, 45 (1992) 288--300.

\litem{[FS]}	Feinstein,~J.F. and Somerset,~D.W.B.,
		Non-regularity for Banach function algebras,
		{\sl Studia Math.}, 141 (2000), 53--68.

\litem{[G]}	Gamelin,~T.W.,
		{\sl Uniform algebras and Jensen measures}.
		London Mathematical Society Lecture Note Series, 32.
		Cambridge University Press, Cambridge-New York, 1978.

\litem{[GL]}	Gamelin, T. W.; Lyons, T. J.,
		Jensen measures for $R(K)$,
		{\it J. London Math. Soc.} (2) 27 (1983), 317--330.

\litem{[H]}	Huber,~A., \"Uber Potentiale, welche auf vorgebenen Mengen verschwinden,
		{\it Comment. Math. Helv.} 43 (1968), 41--50.

\litem{[Ka]}	Karahanjan,~M.I., Some algebraic characterizations of
		the algebra of all continuous functions on a locally connected
		compactum, Math. USSR-Sbornik, 35 (1978), 681--696.

\litem{[Ko]}	K\"orner,~T.W., A cheaper Swiss cheese, {\it Studia Math.},
		83 (1986), 33--36.

\litem{[M]}	McKissick,~R., A nontrivial normal sup norm algebra, {\it Bull.,
		Amer. Math. Soc.} 69 (1963), 391--395.

\litem{[O]}	O'Farrell,~A.G., A regular uniform algebra with a continuous point
		derivation of infinite order, {\sl Bull. London Math. Soc.}, 11 (1979) 41-44.

\litem{[S]}	Stout,~E.L., {\sl The Theory of Uniform Algebras}, Bogden and Quigley, Inc.,
		New York 1971.

\litem{[W]}	Wermer, J.,
		Bounded point derivations on certain Banach algebras.
		{\sl J. Functional Analysis}, 1 (1967) 28--36.
\end